\documentclass[12pt]{article}

\def\triangulo{\hbox{.\kern2pt.\kern-5pt\raise 4pt\hbox{.}}\kern 5pt}

\def\bbb#1{\hbox {{\gordas #1}}}

\font\gordas = msbm10 at 12pt
\def\bbb#1{\hbox {{\gordas #1}}}
\def\a{{\bbb A}}

\def\aa{{\bbb A}}
\def\erre{{\bbb R}}
\def\ce{{\bbb C}}
\def\cee{{\bbb C}}
\def\o{{\bbb O}}
\def\oo{{\bbb O}}

\def\ache{{\bbb H}}

\begin{document}
\noindent
\begin{center}
{\large\bf
Constructing zero divisors in the higher dimensional Cayley-Dickson algebras}\\[.5cm]
Guillermo Moreno
\end{center}
\vglue .5cm
\noindent
{\bf Abstract:} In this paper we give methods to construct zero divisors in the Cayley--Dickson algebras
$\a_n$=$\erre^{2^n}$ for $n$ larger than 4. Also we relate the set of zero divisors with suitable Stiefel manifolds.
\newpage
\noindent
{\bf Introduction:}

\vglue.5cm
In this paper we return to the subject of 97' authors paper [M1] about the description of the zero divisors in the Cayley-Dickson algebras over the real numbers.

This paper must be seen as a sequel of  the papers [M1] and [M2] and related to papers [M3] and [M4].

As we know, $\a_n=\erre^{2^n}$ denotes the Cayley Dickson (C-D)-algebras over the real numbers. For $n=0,1,2,3$ are known as the clasical C-D algebra which correspond to $\a_0=\erre,\a_1=\ce,\a_2=\ache$ and $\a_3=\o$, Real, Complex, Quaternions and Octonions number and because they are normed and alternative algebras they lack of zero divisors.[Sch2] [K-Y].

For $n\geq 4, \a_n$ is no normed, non-alternative (but flexible) algebra and has zero divisors i.e.;    there are nonzero elements  $x$ and $y$ such that $xy=0$.

By definition a zero divisor is a nonzero element $a$ such that there exists nonzero $b$ with $ab=0.$

In [M1]  the zero divisors in $\a_4$ and the zero divisors in $\a_{n+1}=\a_n\times\a_n$ with alternative coordinates are described.

For $n=4$  the set of zero divisors in $\a_4$ of fixed norm can be identified with $V_{7,2}$ the real Stiefel Manifold of two frames in $\erre^7$  and the singular set of $(x,y)$ with $xy=0$ and $||x||=||y||=1$ is homoemorphic to
$G_2$ the exceptional simple Lie group of rank 2.See also [K-Y].

The description of the zero divisors in $\a_4$ is given by the known fibration $$G_2\stackrel{\pi}{\rightarrow}V_{7,2}$$

with fiber $S^3$ since all the nontrivial  annihilators are 4 dimensional.
\vglue .25cm

For $n\geq 5$ there is NO  analogous description. We will show that the zero divisors are in $\a_{n+1}$
and $V_{2^n -1,2}$ are related, but they are not equal and the corresponding singular set has (unknown)
complicated description. See [M3].

\vglue .25cm

In $\S$1 we recall the basic notation and theorems already proved in [M1], [M2] which are necesary for further results.

In $\S$2 we study the basic facts of the linear operators  left and right mutiplication by a fixed (pure) element, as in [M2], and then, we define a suitable $O(2)$-action on the double pure elements of $\a_n.$ 

In $\S$3 we define the {\it Spectrum}  of a non-zero double pure element looking at the structure of the linear operators defined by left and right multiplication by the element.

By definition, the Spectrum  is a set of $(2^{n-2} -1 )$ non negative real numbers attached to each nonzero
double pure element in $\a_n$ for $n\geq 3$ and the presence of zero in the spectrum determines that the element is a zero divisor.

As consequence we see how big can be the annihilator of an element; this result complements very nicely with the recent results in  [B-D-I].

Also we look at the $O(2)$-action on the doubly pure elements showing that: the elements in the same $O(2)$-orbit have the same spectrum and consequently the annihilators of the elements in the $O(2)$-orbit are all equal.

In $\S$4 we study the zero divisors and we give new methods to construct  many of them, noticing that for $n$ larger than 4 the problem becomes extremely difficult.

Also we prove that:All non-zero pure element in $\a_n$ is a component of a zero divisor in $\a_{n+1}$

In $\S$5 we study the relationship between the Stiefel manifold $V_{2^n-1,2}$ [J] and the set of zero divisors in $\a_{n+1}$ noticing that for $n$ larger than 4 they are realted but are "very far" to be equal.

Some of the results of this paper were partially done in a personal manuscript by the author of 2001,(which had restricted circulation) and some also can be  found in [Ch].

We emphasize the initial input and permanent support to study this subject  given by Professor Fred Cohen with whom we are very grateful.
\newpage
\noindent
{\bf 1. Basic definitions and lemmas.}
\vglue.5cm
Here $e_0 =(e_0, 0)$ is the algebra unit in $\aa_{n-1}\times \aa_{n-1}=\aa_n$.
\vglue .25cm
\noindent
{\bf Definition:} $e_{2^{n-1}}=(0, e_0)\in \aa_{n-1}\times \aa_{n-1}=\aa_n$ is {\it the sympletic unit} and we denote it by $\widetilde{e}_0=(0,e_0)$. 
For
$$
a=(a_1, a_2)\in \aa_{n-1}\times \aa_{n-1}=\aa_n\\
$$
$$
\widetilde{a}:=(-a_2, a_1)=(a_1, a_2)(0, e_0)=a\widetilde{e}_0.
$$
is the complexification of $a$.

Because the right multiplication by $\widetilde{e}_0$.   
$$
R_{\widetilde{e}_0}:\aa_n\rightarrow \aa_n
$$
has the matrix $\biggl(\begin{array}{ll}0 &-I \\I &0\end{array}\biggr)$ in $M_{2^{n-1}}(\erre)$ in the canonical basis.

Notice that $\widetilde{\widetilde{a}}=-a$ for all $a\in \aa_n$; {\it the trace} on $\aa_n$ is the linear transformation
$$
t_n :\aa_n\rightarrow \aa_0=\erre
$$
given by $t_n(a)=a+\overline{a}=2$ (Real part of $a$). For $a=(a_1, a_2)\in \aa_{n-1}\times \aa_{n-1}=\aa_n$
$$t_n(a)=t_{n-1}(a_1).$$
\vglue .5cm
\noindent
{\bf Definition:} An element $a$ {\it is pure} in $\aa_n$ if 
$$
t_n(a)=t_{n-1}(a_1)=0
$$
Notice that $a\in \aa_n$ is pure if and only if $a\perp e_0$ and $a\perp b$ if and only if $ab=-ba$ for $a$ and $b$ pure elements in $\aa_n$.
\vglue .5cm
\noindent
{\bf Definition:} An element $a=(a_1, a_2)\in \aa_{n-1}\times \aa_{n-1}=\aa_n$ is {\it doubly pure} if it is pure i.e. $t_n(a)=0$ and also $t_n(\widetilde{a})=-t_{n-1}(a_2)=0$.

Notice that $a$ is doubly pure if and only if $a\in \{e_0, \widetilde{e}_0\}^\perp$.
\vglue .5cm
\noindent
{\bf Notation:} $_0\aa_n =\erre^{2^n-1}=\{e_0\}^\perp$ pure elements in $\aa_n$. $\widetilde{\aa}_n=\erre^{2^n-2}=\{e_0, \widetilde{e}_0\}^\perp$ doubly pure elements in $\aa_n$.
\vglue .5cm
\noindent
{\bf Lemma 1.1.} For $a$ and $b$ in $\widetilde{\aa_n}$ we have that 
\begin{itemize}
\item[1)] $a\widetilde{e}_0=\widetilde{a}$ and $\widetilde{e}_0a=-\widetilde{a}$
\item[2)] $a\widetilde{a}=-||a||^2\widetilde{e}_0$ and $\widetilde{a} a=||a||^2\widetilde{e}_0$ so $a\perp \widetilde{a}$
\item[3)] $\widetilde{a}b=-\widetilde{ab}$ with $a$ pure element.
\item[4)] $a\perp b$ if and only if $\widetilde{a}b+\widetilde{b}a=0$
\item[5)] $\widetilde{a}\perp b$ if and only if $ab=\widetilde{b}\widetilde{a}$
\item[6)] $\widetilde{a}b=a\widetilde{b}$ if and only if $a\perp b$ and $\widetilde{a}\perp b$.
\end{itemize}
\vglue .25cm
\noindent
{\bf Proof:} Notice that $a$ is pure if $\overline{a}=-a$ and if $a=(a_1, a_2)$ is doubly pure then $\overline{a}_1=-a_1$ and $\overline{a}_2=-a_2$.
\begin{itemize}
\item[1)] $\widetilde{e}_0a=(0, e_0)(a_1, a_2)=(-\overline{a}_2, \overline{a}_1)=(a_2, -a_1)=-(-a_2,a_1)=-\widetilde{a}$.
\item[2)] $a\widetilde{a}=(a_1, a_2) (-a_2, a_1)=(-a_1a_2+a_1a_2, a_1^2+a^2_2)=-||a||^2\widetilde{e}_0$

$\widetilde{a}a=(-a_2, a_1)(a_1, a_2)=(-a_2a_1+a_2a_1, -a^2_2-a^2_1)=||a||^2\widetilde{e}_0$.Now since $-2\langle \widetilde{a}, a\rangle =a\widetilde{a}+\widetilde{a}a=0$ we have $a\perp \widetilde{a}$.
\item[3)] $\widetilde{a}b=(-a_2, a_1)(b_1, b_2)=(-a_2b_1+b_2a_1, -b_2a_2 -a_1b_1)$. So $\widetilde{\widetilde{a}b}= (a_1b_1+b_2a_2, b_2a_1-a_2b_1)=(a_1, a_2)(b_1, b_2)=ab$ then $-\widetilde{a}b=\widetilde{\widetilde{\widetilde{a}b}}=\widetilde{ab}$.

Notice that in this proof we only use that $\overline{a}_1=-a_1$ i.e $a$ is pure and $b$ doubly pure.
\item[4)] $a\perp b\Leftrightarrow ab+ba=0\Leftrightarrow ab=-ba\Leftrightarrow\widetilde{ab}=-\widetilde{ba}$

$\Leftrightarrow -\widetilde{a}b=\widetilde{b}a\Leftrightarrow \widetilde{a}b+\widetilde{b}a=0$ by (3).
\item[5)] $\widetilde{a}\perp b \Leftrightarrow \widetilde{\widetilde{a}}b+\widetilde{b}\widetilde{a}=0$ (by (4)) $\Leftrightarrow -ab+\widetilde{b}\widetilde{a}=0$.
\item[6)] If $\widetilde{a}\perp b$ and $a\perp b$ then
$$
\widetilde{a}b=-\widetilde{ab}=\widetilde{ba}=-\widetilde{b}a=a\widetilde{b}.
$$
Conversly, put $a=(a_1, a_2)$ and $b=(b_1, b_2)$ in $\aa_{n-1}\times \aa_{n-1}$ and define $c:=(a_1b_1+b_2a_2)$ and $d:=(b_2a_1-a_2b_1)$ in $\aa_{n-1}$.

Then $a\widetilde{b}=(a_1, a_2)(-b_2, b_1)=(-a_1b_2+b_1a_2, b_1a_1+a_2b_2)$ so $a\widetilde{b}=(-\overline{d}, \overline{c})$.
\end{itemize}
Now $ab=(a_1, a_2)(b_1, b_2)=(a_1b_1+b_2a_2, b_2a_1-a_2b_1)=(c,d)$ then $\widetilde{ab}=(-d,c)$ so $\widetilde{a}b=(d, -c)$. Suppose that $a\widetilde{b}=\widetilde{a}b$ so $\overline{c}=-c$ and $d=-\overline{d}$ then 
\begin{eqnarray*} 
t_n(ab)&=&t_{n-1}(c)=c+\overline{c}=0\;\;\hbox{\rm and}\;\;a\perp b\\
t_n(\widetilde{a}b)&=&t_{n-1}(d)=d+\overline{d}=0\;\;\hbox{\rm and}\;\;\widetilde{a}\perp b.
\end{eqnarray*}

\hfill Q.E.D.
\vglue .5cm
\noindent
{\bf Corollary 1.2} For each $a\neq 0$ in $\widetilde{\aa}_n$. The fourth dimensional vector subspace generated by $\{e_0, \widetilde{a}, a, \widetilde{e}_0\}$ is a copy of $\aa_2 =\ache$ the Hamilton quaternions. We denote it by $\ache_a$.
\vglue .25cm
\noindent
{\bf Proof:} We suppose that $||a||=1$, otherwise we take ${a\over ||a||}$. Construct the following multiplication table
\vglue .5cm
\begin{center}
\begin{tabular}{c|cccc}
 & $e_0$ & $\widetilde{a}$ & $a$ & $\widetilde{e}_0$\\ \hline
$e_0$ & $e_0$ & $\widetilde{a}$ & $a$ & $\widetilde{e}_0$ \\
$\widetilde{a}$ & $\widetilde{a}$ & $-e_0$ & $+\widetilde{e}_0$ & $-a$\\
$a$ & $a$ & $-\widetilde{e}_0$ & $-e_0$ & $\widetilde{a}$\\
$\widetilde{e}_0$ & $a$ & $\widetilde{e}_0$ & $-\widetilde{a}$ & $-e_0$
\end{tabular}
\end{center}
By lemma 1.1. $a\widetilde{e}_0=\widetilde{a}$; $\widetilde{e}_0a=-\widetilde{a}$; $\widetilde{a}\widetilde{e}_0=\widetilde{\widetilde{a}}=-a$; $\widetilde{e}_0 \widetilde{a}=-\widetilde{\widetilde{a}}=a$; $a\widetilde{a}=-\widetilde{e}_0$ and $\widetilde{a}a=\widetilde{e}_0$.

This multiplication table is the one of the Hamilton quaternions identifying $e_0\leftrightarrow 1$, $\widetilde{a}\leftrightarrow \hat{i}$, $a\leftrightarrow \hat{j}$ and $\widetilde{e}_0\leftrightarrow \hat{k}$.

\hfill Q.E.D.

 Artin theorem [Sch1] said that: ``Any two elements in $\aa_3=\oo$ generate an associative subalgebra''.

Corollary 1.2 is an analogous to this for $\aa_n$ with $n\geq 4$.
\vglue .5cm
\newpage
\noindent
{\bf 2. Left and Right multiplication.}
\vglue1cm
For $x$ and $y$ in $\aa_n\;\;\langle x, y\rangle$ denotes the standard Euclidean inner product that in terms of the Cayley-Dickson multiplication is given by 
$$
t_n(x \overline{y})=x\overline{y}+y\overline{x}=2\langle x, y\rangle .
$$
In [1] and [8 ] is proved that
\begin{eqnarray*}
\langle a x, b\rangle & =& \langle x, \overline{a}b\rangle\qquad\hbox{\rm and}\\
\langle a, xb \rangle &=& \langle a\overline{b}, x\rangle
\end{eqnarray*}
for all $a,b$ and $x$ in $\aa_n$.

Therefore if $L_a$ and $ R_b :\aa_n \rightarrow \aa_n$ denotes the linear transformation left and right multiplication by $a$ and $b$ respectively then $L^T_a=L_{\overline{a}}$ and $R^T_b=R_{\overline{b}}$.

So for $a$ and $b$ pure elements $L^T_a=-L_a$ and $R^T_b=-R_b$ i.e. $L_a$ and $R_b$ are skew--symmetric linear transformations then $(-L^2_a)$ and $(-R^2_b)$ are symmetric definite nonegatives linear transformations.
\vglue .5cm
\noindent
{\bf Lemma 2.1} For $a$  a doubly pure element in $\aa_n,\;\;n\geq 3$
$$
R_a R_{\widetilde{e}_0}+R_{\widetilde{e}_0}R_a=0\;\;\;\hbox{\rm and}\;\;\;L_aL_{\widetilde{e}_0}+L_{\widetilde{e}_0}L_a=0.
$$
{\bf Proof:} Using Lemma 1.1 (3).

$(R_aR_{\widetilde{e}_0}+R_{\widetilde{e}_0}R_a)(x)=(x\widetilde{e}_0)a+(xa)\widetilde{e}_0=\widetilde{x}a+\widetilde{xa}=0$ for $x\in\,_0\aa_n$. If $x=re_0$ for $r\in \erre\;\;\widetilde{x}a+\widetilde{xa}=r\widetilde{e}_0a+\widetilde{(re_0)a}=r(-\widetilde{a}+\widetilde{a})=0$. Then $(R_aR_{\widetilde{e}_0}+R_{\widetilde{e}_0}R_a)(x)=0$ for all $x\in \aa_n$.

On the other hand we have that for $y=(y_1, y_2)$
$$
\widetilde{e}_0y=(0, e_0)(y_1, y_2)=(-\overline{y}_2, \overline{y}_1)
$$
then 
\begin{eqnarray*}
\widetilde{e}_0 (ax) &=& (0, e_0)(a_1 x_1-\overline{x}_2 a_2, x_2a_1+a_2\overline{x}_1)\\
&=&(-\overline{(x_2a_1+a_2\overline{x}_1)}, \overline{(a_1x_1-\overline{x}_2a_2)})\\
&=& (-(\overline{a}_1\overline{x}_2+x_1\overline{a}_2), (\overline{x}_1 \overline{a}_1-\overline{a}_2x_2))\;\;\;\hbox{\rm and  because}\;a\in \widetilde{\aa_n}\\
&=&(a_1\overline{x}_2+x_1a_2, -\overline{x}_1a_1+a_2x_2)\;\;\;\hbox{\rm for}\;\;x=(x_1,x_2) \in \aa_n
\end{eqnarray*}
and $a(\widetilde{e}_0x)=(a_1, a_2)(-\overline{x}_2, \overline{x}_1)=(-a_1\overline{x}_2-x_1a_2, \overline{x}_1a_1-a_2\overline{x}_2)$ 
 therefore $(L_{\widetilde{e}_0}L_a+L_aL_{\widetilde{e}_0})(x)=\widetilde{e}_0 (ax)+a(\widetilde{e}_0 x)=0$ for $x\in \aa_n$

\hfill Q.E.D.
\vglue .5cm
\noindent
{\bf Theorem 2.2.} For $a\in \aa_n$ pure element $\;L^2_a=R^2_a$.
\vglue .25cm
\noindent
{\bf Proof:} (Case doubly pure). For $a\in \widetilde{\aa}_n$, $\aa_n=\ache_a\oplus \ache^\perp _a$ so if $x\in \ache_a$ then (because $\ache_a$ is associative)$$
a(ax)=a^2x=-||a||^2x=-x||a||^2=xa^2=(xa)a.
$$
If $x\in \ache^\perp $ then $a\perp x$ and by flexibility $a(ax)=-a(xa)=-(ax)a=(xa)a$ so $L^2_a=R^2_a$. 

(General case).
Suppose that $a\in \widetilde{\aa}_n$ and $r\in \erre$ so 
\begin{eqnarray*}
L^2_{a+r\widetilde{e}_0}&=&(L_a+rL_{\widetilde{e}_0})^2=L^2_a+r^2L^2_{\widetilde{2}_0}+r(L_aL_{\widetilde{e}_0}+L_{\widetilde{e}_0}L_a)\\
R^2_{a+r\widetilde{e}_0}&=&(R_a+rR_{\widetilde{e}_0})^2=R^2_a+r^2R^2_{\widetilde{e}_0}+r(R_aR_{\widetilde{e}_0}+R_{\widetilde{e}_0}R_a).
\end{eqnarray*}
By the doubly pure case $L^2_a=R^2_a$, also $L^2_{\widetilde{e}_0}=R^2_{\widetilde{e}_0}=-I$ by lemma 2.1 $L_a L_{\widetilde{e}_0}+L_{\widetilde{e}_0}L_a=R_aR_{\widetilde{e}_0}+
R_{\widetilde{e}_0} R_a=0$ therefore $R^2_{a+r\widetilde{e}_0}=L^2_{a+r\widetilde{e}_0}=L^2_a-r^2I=R^2_a-r^2I$.

\hfill Q.E.D.
\vglue .5cm
\noindent
{\bf Remark:} For $0\neq a$ pure $L_a$ and $R_a$ are quite different because the center of $\aa_n$ is $\erre{e_0}$.
\vglue .25cm
\noindent
{\bf Notation:} For $a$ and $b$ non--zero pure elements.
$$
{\cal A}=L^2_a+R^2_b \;\;\;\hbox{\rm and}\;\;\;S=(a, -, b)=R_bL_a-L_aR_b.
$$
are linear transformations ${\cal A}, S: \aa_n \rightarrow \aa_n$.

Notice that ${\cal A}$ is the sum of two symmetric non--positive definite linear transformations so ${\cal A}$ is also symmetric non--postive definite i.e. $L^2_a\leq 0$ and $R^2_b\leq 0$ implies that ${\cal A}=L^2_a+R^2_b\leq 0$ and $S=[R_b, L_a]$ is  skew--symmetric.
\vglue .5cm
\noindent
{\bf Theorem 2.3} For $a$ and $b$ (non--zero) pure elements in $\aa_n$
$$ 
L^2_{(a,b)}:\aa_{n+1}\rightarrow \aa_{n+1} 
$$
is given by
$$
L^2_{(a,b)} =\left(
\begin{array}{ll}
{\cal A} & - S\\
S & {\cal A}
\end{array}\right)
$$
i.e. $L^2_{(a,b)}(x,y)=({\cal A}(x)-S(y), {\cal A}(y)+S(x))$.
\vglue .25cm
\noindent
Proof: (By direct calculation).
\begin{eqnarray*}
(a,b)[(a,b)(x,y)] &=& (a,b)(ax-\overline{y}b, ya+b\overline{x})\\
&=& (a(ax)-a(\overline{y}b)-\overline{(ya+b\overline{x})}b, (ya+b\overline{x})a+b\overline{(ax-\overline{y}b)}\\
&=&(a(ax)-a(\overline{y}b)+(a\overline{y})b-(x\overline{b})b, (ya)a+(b\overline{x})a-b(\overline{x}a)-b(\overline{b}y))\\
&=&(L^2_a(x)+(a, \overline{y}, b)+R^2_b(x), R^2_a (y)+(b, \overline{x}, a)+L^2_b (y))\\
&=&({\cal A}(x)-(a,y,b), {\cal A}(y)+(a, x, b))\\
&=&({\cal A}(x)-S(y), {\cal A}(y)+S(x)).
\end{eqnarray*}

\hfill Q.E.D.
\noindent
{\bf Corollary 2.4} For $r$ and $s$ real numbers with $r^2+s^2=1$ and $a$ and $b$ pure elements in $\aa_n$ we have that 
\begin{enumerate}
\item[1)] $L^2_{(ra-sb, sa+rb)}=L^2_{(a,b)}$ in particular $L^2_{(-b,a)}=L^2_{(a,b)}$
\item[2)] $L^2_{(ra+sb, sa-rb)}=L^2_{(b,a)}$ in particular $L^2_{(a, -b)}=L^2_{(b,a)}$
\end{enumerate}
\noindent
{\bf Proof:} 1)
\begin{eqnarray*}
L^2_{ra-sb}&=&(rL_a-sL_b)^2 =r^2L^2_a+s^2L^2_b-rs (L_aL_b+L_b L_a)\\
L^2_{sa+rb}&=&(sL_a+rL_b)^2 =s^2L^2_a+r^2L^2_b+sr(L_aL_b+L_bL_a).
\end{eqnarray*}
Thus $L^2_{ra-sb}+L^2_{sa+rb}=(r^2+s^2)L^2_a+(r^2+s^2)L^2_b=L^2_a+R^2_b:={\cal A}$.

On the other hand for all $x\in \aa_n$
\begin{eqnarray*}
(ra-sb, x, sa+rb)&=&(ra, x, sa)+(ra, x,rb)-(sb, x, sa)-(sb, x, rb)\\
&=&0+r^2(a,x,b)-s^2(b,x,a)-0\\
&=&(r^2+s^2)(a,x,b)\\
&:=& S(x)
\end{eqnarray*}
by Theorem 2.2 we are done with 1).

To show 2) we observe that ${\cal A}=L^2_a +R^2_b=L^2_b+R^2_a$ and 
$$
S=(a, -,b)=-(b, -a)\;\;\hbox{\rm  so}\;\; L^2_{(b,a)}=\left(
\begin{array}{ll}
{\cal A} & S\\
-S & {\cal A}
\end{array}\right)
$$
and the prove is analogous to case 1).

\hfill Q.E.D.
\vglue .5cm
\noindent
{\bf Corollary 2.5} For $\alpha =(a,b)\in \aa_n\times \aa_n=\aa_{n+1}$ doubly pure we have that:
$$
L_\alpha L_{\widetilde{\alpha}}+L_{\widetilde{\alpha}}L_\alpha =0
$$
and
$$
R_\alpha R_{\widetilde{\alpha}}+R_{\widetilde{\alpha}}R_\alpha =0.
$$
\noindent
{\bf Proof:} Suppose that $r$ and $s$ are in $\erre$ and $r^2+s^2=1$ so $r\alpha +s\widetilde{\alpha}=r(a,b)+s(-b,a)=(ra-sb, rb+sa)$ then $L^2_{r\alpha +s\widetilde{\alpha}}=(r^2+s^2)L^2_\alpha$ by Corollary 2.4.

On the other hand 
\begin{eqnarray*}
L^2_{r\alpha +s\widetilde{\alpha}} &=&(r L_\alpha +sL_{\widetilde{\alpha}})^2 = r^2L^2_\alpha +s^2L^2_{\widetilde{\alpha}}+rs (L_\alpha L_{\widetilde{\alpha}}+L_{\widetilde{\alpha}}L_\alpha)\\
&=& (r^2+s^2)L^2_\alpha \;\;\;\;\hbox{\rm because}\;\;\;L^2_\alpha =L^2_{\widetilde{\alpha}}.
\end{eqnarray*}
so $L_\alpha L_{\widetilde{\alpha}}+L_{\widetilde{\alpha}}L_\alpha =0$.

Similarly $R_\alpha R_{\widetilde{\alpha}}+R_{\widetilde{\alpha}}R_\alpha =0$.

\hfill Q.E.D.
\vglue .5cm

{\bf Remark} Notice that Lemma 2.1 and Corollary 2.5 give us a way to define {\it Cuaternionic Structures} 
on $\a_n$ for n bigher than 2 via alternative elements of norm one, because for such an element $a$
 we have also that $L^2_a=R^2_a= -I.$ (See [Po]). 
 \vglue .5cm

Now the theory of Schur complement for partitioned matrices ([ Z ] Chapter 2 and [H-J]) tell us that.

If ${\cal A}$ is an invertible matrix then
$$
\hbox{\rm det}(L^2_{(a,b)})=\;\hbox{\rm det}\;({\cal A})\,\hbox{\rm det}\,({\cal A}+S{\cal A}^{-1}S)
$$
On the other hand ${\cal A}$ is symmetric definite non--positive i.e. ${\cal A}\leq 0$ because also $L^2_a\leq 0$ and $R^2_b\leq 0$ then $-{\cal A}\geq 0$, $-L^2_a\geq 0$ and $-R^2_b\geq 0$ and det$(-{\cal A})=$det$(-L^2_a-R^2_b)\geq$det$(-L^2_a)+$det$(-R^2_b)$ but ${\cal A}$, $L^2_a$ and $R^2_b$ are matrices of order $2^n$ so det$({\cal A})\geq$ det$(L_a)^2+$ det$(R_b)^2$.

Therefore if det$({\cal A})=0$ then det$(L_a)=0$ and det$(R_b)=0$ and ${\cal A}$ is invertible if, either, $L_a$ or $R_b$ is invertible. 

Suposse now that $-{\cal A}>0$. i.e. $L_a$ or $R_b$ are invertible then $-{\cal A}^{-1}>0$ and since $S^T=-S$. $(S{\cal A}^{-1}S)^T=S^T({\cal A}^{-1})^TS^T=S{\cal A}^{-1}S$ so $S{\cal A}^{-1}S$ is symmetric and $\langle S{\cal A}^{-1}S(x), x\rangle =\langle -{\cal A}^{-1}S(x), S(x)\rangle \geq 0$ for all $x\in \aa_n$, because $-{\cal A}^{-1}>0$, so $(S{\cal A}^{-1}S)\geq 0$ and $(S{\cal A}^{-1}S)(x)=0$ if and only if $S(x)=0$. 

We resume this discussion in 
\vglue .5cm
\noindent
{\bf Theorem 2.6} For $(a,b)\in \,_0\aa_n\times \,_0\aa_n=\widetilde{\aa}_{n+1}$ for $n\geq 3$. If $L_a$ or $R_b$ are invertible and $\lambda =-1$ is no an eigenvalue of $({\cal A}^{-1}S)^2$ then $L_{(a,b)}$ is invertible in $\aa_{n+1}$.
\vglue .25cm
\noindent
{\bf Proof:} We already show that $L_a$ or $R_b$ invertible implies that ${\cal A}$ is invertible.

Now $L_{(a,b)}$ is invertible if and only if $L^2_{(a,b)}$ is invertible and, by Schur complement, that happen if and only if det$({\cal A}+S{\cal A}^{-1}S)\neq 0$ but det$({\cal A}+S{\cal A}^{-1}S)=$ det$({\cal A}^{-1})$ det$(I+({\cal A}^{-1}S)^2)$ so if $\lambda =-1$ is an eigenvalue of $({\cal A}^{-1}S)^2$ then det$(I+({\cal A}^{-1}S)^2)=0$.

\hfill Q.E.D.
\vglue .5cm
This theorem is useful when ${\cal A}$ has a simple expression for instance when $a$ and $b$ are alternative elements so ${\cal A}=L^2_a+R^2_b=(a^2+b^2)I$ so ${\cal A}^{-1}=(a^2+b^2)^{-1}I$ and ${\cal A}^{-1}S=(a^2+b^2)^{-1}S$ then $L^2_{(a,b)}$ is singular if there exists $x\neq 0$ in $\aa_n$ such that $S^2(x)=-(a^2+b^2)^2x$. (See [8] $\S$ 2).
\newpage
\noindent
{\Large\bf $\S$ 3. The spectrum of a doubly pure element.}

\vglue .5cm
Let $a\in \aa_n$ be a doubly pure element with $||a||\neq 0$

By the diagonalization theorems of skew--symmetric matrices [Pr] we know that there exists an orthogonal basis with respect to which $L_a$ has the form
$${\rm diag}
\left(
\wedge_1,\wedge_2, \cdot ,\cdot, \cdot,\wedge_k,0\cdot,\cdot,\cdot, 0
\right)
$$
where $\wedge_i =\left(\begin{array}{cc}0 &-\lambda_i\\\lambda_i &0
\end{array}\right)$ for $\lambda_i\geq 0$ in $\erre$. So with respect to the same basis $L^2_a$ has the form
$${\rm diag}
\left(
-\lambda^2_1, -\lambda^2_1, \cdot, \cdot, \cdot ,-\lambda^2_k, -\lambda ^2_k,
0, 0\cdot, \cdot, 0, 0\right)
$$
Now, recall that  $L_a(\ache_a)=\ache_a$ and $L_a(\ache^\perp _a)\subset \ache^\perp_a$.
\vglue .25cm

Since $\ache_a$ is associative we have that
$$
L^2_a|_{\ache _a}=a^2I_{4\times 4}
$$
Therefore the first two $\lambda 's$,(i.e. $\lambda_1$ and$\lambda_2$) are equal to $||a||$ so we restrict ourselves to
$$
L_a, L^2_a:\ache^\perp_a\rightarrow \ache^\perp_a.
$$
Define for $\lambda\geq 0$ and $a$ doubly pure with $||a||=1$
$$
V_\lambda =\{x\in \ache^\perp_a|a(ax)=-\lambda^2 x\}.
$$
\noindent
{\bf Theorem 3.1} For $0\neq x\in V_\lambda$ the set $\{x, (-1/\lambda) (ax), (-1/\lambda) ( \widetilde{ax}), \widetilde{x}\}$ is an orthogonal set in $V_\lambda$ for nonzero $\lambda$ and the dimension of $V_\lambda$ is congruent with $0$ mod $4$.
\vglue .25cm
\noindent
{\bf Proof:} If $0\neq x\in V_\lambda$ then $\widetilde{x}\in V_\lambda$ and $\widetilde{a(ax)}=-\lambda^2 \widetilde{x}$. 
Also 
$$
a(a((-1/\lambda) (ax)))=(-1/\lambda ) a(a(ax))=(-1/\lambda) (-\lambda^2ax)=-\lambda^2(-1/\lambda )( ax)
$$
and $(-1/\lambda) (ax)\in V_\lambda$ therefore $(-1/\lambda) (\widetilde{ax})\in V_\lambda$ and by construction \\
$\{x, (-1/\lambda) (ax), (-1/\lambda) (\widetilde{ax}), \widetilde{x}\}$ is an orthonormal set. Now take $0\neq y\in (\ache_a\oplus \{x,(1/\lambda) (ax), (-1/\lambda) (\widetilde{ax}), \widetilde{x}\})^\perp$ 

Therefore $y$ is doubly pure and  orthogonal to $a, \widetilde{a}, x, ax,\widetilde{ax}$ and $\widetilde{x}$ using that $L^T_a=-L_a$ we have that:
\begin{eqnarray*}
\langle ay, x\rangle &=&\langle y,- ax\rangle =0\;\;\;\hbox{\rm so}\;\;\;ay\perp x.\\
\langle ay,ax\rangle &=&\langle y, -a(ax)\rangle =\langle y, -(-\lambda^2 x)\rangle =\lambda^2\langle y,x\rangle =0\\
\langle ay,\widetilde{ax}\rangle &=&\langle ay, -\widetilde{a}x\rangle =\langle \widetilde{a}(ay), x\rangle =\langle -\widetilde{a(ay)}, x\rangle\\
&=&\langle a(ay), \widetilde{x}\rangle =\langle y, a(a\widetilde{x})\rangle =\langle y, -\lambda^2\widetilde{x}\rangle =0\\
\langle ay, \widetilde{x}\rangle &=&\langle -\widetilde{ay}, x\rangle =\langle \widetilde{a}y, x\rangle =\langle y, -\widetilde{a}x\rangle =\langle y, \widetilde{ax}\rangle =0\\
\hbox{\rm Similarly}&&\langle \widetilde{y}, x\rangle =\langle -y,\widetilde{x}\rangle=0; \langle \widetilde{y}, ax\rangle =\langle y, -\widetilde{ax}\rangle=0\\
\langle \widetilde{y}, \widetilde{ax}\rangle &=&\langle y, ax\rangle =0
\end{eqnarray*}
Finally $\;\;\langle \widetilde{ay}, x\rangle =\langle -\widetilde{a}y, x\rangle =\langle y, \widetilde{a}x\rangle =0$ and 
\begin{eqnarray*}
\langle \widetilde{ay}, ax\rangle &=&\langle -\widetilde{a}y, ax\rangle =\langle y, \widetilde{a}(ax)\rangle =\langle y, -\widetilde{a(ax)}\rangle =\\
&=&\langle y, +\lambda^2\widetilde{x}\rangle =\lambda^2\langle y, \widetilde{x}\rangle =0\\
\langle \widetilde{ay}, \widetilde{ax}\rangle &=&\langle ay, ax\rangle =0
\end{eqnarray*}
Therefore for any $0\neq y\in \ache^\perp_a$ with $y\in \{x, (-1/\lambda) ax,( 1/\lambda)x, \widetilde{ax}, \widetilde{x} \}^{\perp }$ we have that $Span\{y, \widetilde{y} , ay, \widetilde{ay}\}\cap Span\{x, (-1/\lambda) ax, (-1/\lambda) \widetilde{ax}, \widetilde{x}\}=\{0\}$.

\hfill Q.E.D.
\vglue .5cm

Therefore the eigenvalues of $L^2_a:\ache _a^\perp \rightarrow \ache^\perp _a$ have multiplicity congruent with $0$ mod $4$ and 
$$
\ache^\perp_a=V_{\lambda_1} \oplus V_{\lambda_2}\oplus \cdots \oplus V_{\lambda_k}\qquad k={2^n-4\over 4}=2^{n-2}-1
$$
\vglue .25cm
\noindent
{\bf Definition:} The {\it Spectrum} of $a$ a non--zero doubly pure element, Spec$(a)$, is the set of $(2^{n-2}-1)$ eigenvalues (of multiplicity congruent with $0$ mod $4$) of 
$$-L^2_{a\over ||a||}:\ache^\perp_a\rightarrow \ache^\perp_a$$.

By definition the elements in Spectra are nonegative real numbers.
\vglue .25cm

For instance if  $a\neq 0$ in $\widetilde{\aa_3}$ then
$$
\hbox{\rm Spec}\;(a)=\{1\}.
$$
Recall that $a\in \aa_n$ is alternative if $a(ax)=a^2x$ for all $x$ so for $0\neq a\in \aa_n$ alternative.
$$
\hbox{\rm Spec}\;(a)=\{1, 1, \ldots , 1\}.
$$
Also $L_a$ is non--singular if and only if $0\not\in$ Spec$(a)$.
\vglue .25cm
\noindent
{\bf Example:} Take $e_1$ and $e_2$ basic elements in $\aa_3$.

Define $a=(e_1, e_2)\in \aa_4\;\;$ so $||a||=\sqrt{2}$ and
\begin{eqnarray*}
(e_1, e_2)(-e_4, e_7)&=&(-e_1e_4+e_7e_2, e_7e_1+e_2e_4)=(-e_5+e_5, -e_6+e_6)\\
&=&(0,0)
\end{eqnarray*}
so $0\in$ Spec$(a)$ and $L_{(e_1, e_2)}$ is non--invertible.

Also
\begin{eqnarray*}
(e_1, e_2)(e_2, e_1)&=&(2e_1 e_2, 0)\qquad\hbox{\rm and}\\
L^2_{(e_1, e_2)}(e_2, e_1)&=&2(e_1 (e_1e_2), e_2 (e_2, e_1))=-2(e_2, e_1).
\end{eqnarray*}
and $1\in$ Spec$(a)$.

The spectrum of $a$ has three elements we already calculate two $\{0,1\}$.

Consider 
\begin{eqnarray*}
(e_1, e_2)(e_7, -e_4)&=&(e_1e_7-e_4e_2, -e_4e_1-e_2e_7)\\
&=&(e_6+e_6, e_5+e_5)\\
&=&(2e_6, 2e_5).
\end{eqnarray*}
so
\begin{eqnarray*}
L^2_{(e_1, e_2)}(e_7, -e_4)&=&2(e_1, e_2)(e_6, e_5)\\
&=&2[(e_1e_6+e_5e_2, e_5e_1-e_2e_6)]\\
&=&2[(-e_7-e_7, e_4+e_4)]\\
&=&2[(-2e_7. 2e_4)]\\
&=&-4(e_7, -e_4).
\end{eqnarray*}
Since $|(e_1, e_2)|=\sqrt{2}$ then $L^2_{a\over ||a||}=({1\over \sqrt{2}})^2L^2_{(e_1, e_2)}$ so 
$$
-{1\over 2} L^{2}_{(e_1, e_2)}(e_7,-e_4)={1\over 2}(4)(e_7, -e_4)=2(e_7, -e_4)
$$
so  Spec$((e_1,e_2))=\{0,1,2\}$.
\vglue .25cm
\noindent
{\bf Remarks:} We can generalize the argument in last example in the following way:

Suppose that $a$ and $b$ are alternative elements in $_0\!\aa_n\;\;n\geq 3$ such that 
$|a|=|b| \neq 0$ and $a$ orthogonal to $b$ then

(i) $1\in$ Spec$((a,b))$ in $\aa_{n+1}$.
\vglue .25cm
(ii) If also $(a,b)(x,y)=(0,0)$ for some $x$ and $y$ in $_0\!\aa_n$ then $2\in$ Spec$((a,b))$ realized by $(y,x)$ and therefore $\{0,1,2\}\subset$ Spec$((a,b))$.

It seems to us that (ii) is an expresion of a"Mirror Symmetry" sourrended by the zero divisors in $\a_{n+1}$: 

"If $(x,y)$ left  $0$ in the spectrum of $(a,b)$ then $(y,x)$ leaves $2$ in the spectrum of $(a,b)$.

"If $(x,y)$ left $1$ in the spectrum of $(a,b)$ then $(y,x)$ leaves $1$ in the spectrum of $(a,b)$.

This "Mirror Symmetry" is "Broken" in the Singular set defined by The zero set of the Hopf map
(see [M3])
$$\{(x,y)| xy=0 , ||x||=||y||=1\}$$
Because $(x,y)(y,x)=(2xy,||y||^2 -||x||^2)$ for $x$ and $y$ pure non zero elements in $\a_n$,
"The Hopf construction map" in terms of the algebra $\a_{n+1}$.
\vglue .25cm
{\bf Remark:}  For $r$ and $s$ in $\erre$ with $r^2 +s^2 =1$ and $(a,b)$ doubly pure in $\a_{n+1}$
 $$\phi(a,b)=(ra\mp sb,sa\pm rb)$$
 
 This define an $O(2)$-action on $\widetilde{\aa}_{n+1}$ and by Corollary 2.4 we have that
 $$ Spec(\phi(a,b))=Spec((a,b))$$
\vglue .5cm
\noindent
Now we "measure" the lack of the Normed property in  $\a_n$ in terms of Spectra.
\vglue .25cm

{\bf Definition:} An element $0\neq a\in \,_0\!\aa_n$ (i.e. $a$ is pure) {\it is normed with $x\in \aa_n$} if $||a||||x||=||ax||$. We also define that {\it $a$ is normed} if $||ax||=||a||||x||$ for all $x\in \aa_n$.
\vglue .5cm
\noindent
{\bf Theorem 3.2} $0\neq a\in \,_0\!\aa_n$ is normed if and only if $a$ is alternative i.e. $a(ax)=a^2x$ for all $x$ in $\aa_n$. 
\vglue .25cm
\noindent
{\bf Proof:} Suppose that $0\neq a\in \,_0\!\aa_n$ is alternative so 
\begin{eqnarray*}
||a||^2 ||x||^2&=&-a^2\langle x, x\rangle =\langle -a^2x, x\rangle =\langle -a(ax), x\rangle =\langle ax, ax\rangle\\
&=&||ax||^2\qquad\hbox{\rm for all}\;\;\;x\in \aa_n.
\end{eqnarray*}
Conversly suppose that $||ax||=||a||||x||$ for all $x\in \aa_n$ then $\langle -a^2x,x\rangle =\langle ax, ax\rangle =\langle -a(ax), x\rangle$ and $\langle a^2x-a(ax), x\rangle =0$ for all $x\in \aa_n$.

But $(a^2I-L^2_a)$ is a symmetric linear transformation then all its eigenvalues are of the form $\langle (a^2I-L^2_a)(x),x\rangle$ so $a^2I-L^2_a=0$.

\hfill Q.E.D.
\vglue .25cm
\noindent
{\bf Remark:} Notice that the properties for an element of being normed and alternative coincide globally but no locally i.e. we say that $a\in \,_0\aa_n$ {\bf alternate with $x$} for $x\in \aa_n$ if $a(ax)=a^2x$ so $a\in \,_0\aa_n$ is an alternative element if $a$ alternates with any element in $\aa_n$ (see[M2]).

Now if $a$ alternate with $x$ then $a$ is normed with $x:||ax||^2=\langle ax,ax\rangle =\langle -a(ax),x\rangle =\langle -a^2x, x\rangle =||a||^2\cdot||x||$.

But the converse is no necessarily true.
\vglue .25cm
\noindent
{\bf Example:} Take $a=e_1$, $b=e_2$ and $\widetilde{e}_0=e_4$ in $\aa_3$. Define $\alpha= (e_1, e_2)$ and $\varepsilon =(\widetilde{e}_0,0)$ in $\aa_3\times \aa_3=\aa_4$. Since $\varepsilon$ is alternative element in $\aa_4$ we have that $||\alpha \varepsilon||^2=\langle \alpha \varepsilon, \alpha \varepsilon \rangle =\langle -(\alpha \varepsilon)\varepsilon, \alpha\rangle =\langle -\alpha \varepsilon^2, \alpha\rangle =||\alpha ||^2||\varepsilon||^2=||\alpha||^2$ and $\alpha $ is normed with $\varepsilon$.

On the other hand using th.2.3
\begin{eqnarray*}
L^2_\alpha(\varepsilon) &=&((L^2_{e_1}+R^2_{e_2})(\widetilde{e}_0)-0, 0+(e_1, \widetilde{e}_0, e_2))\\
&=&((-2I)(e_4), (e_1 e_4)e_2-e_1(e_4e_2))\\
&=&(-2e_4, e_5e_2-e_1(-e_6))\\
&=&(-2e_4, -e_7-e_7)\\
&=&(-2e_4, -2e_7)\\
&=&-2\varepsilon -2e_{15}\qquad\hbox{\rm in}\;\;\aa_4.
\end{eqnarray*}
But $\alpha^2\varepsilon =-2\varepsilon$ so $(\alpha, \alpha, \varepsilon)=\alpha^2\varepsilon -L^2_\alpha (\varepsilon)=2e_{15}$ and $\alpha$ does not alternate  with $\varepsilon$.

\vglue .5cm
\noindent
{\bf Proposition 3.3} For $a\in \widetilde{\aa}_n$ with $a\neq 0$ and $\lambda \in$ Spec$(a)$. 

Suppose that $0\neq x\in V_\lambda \subset \ache^\perp_a \subset \widetilde{\aa}_n$ then 

(1)$||ax||=\lambda ||a||||x||$.

(2) If $x$ is an alternative element then $\lambda =1$.

\vglue .25cm
\noindent
{\bf Proof:} (1) Assuming that $||a||=1$ we have that
$$a(ax)=-\lambda^2 x.$$
so
$$||ax||^2=\langle ax,ax\rangle =-\langle a(ax),x\rangle =\langle \lambda^2x,x\rangle =\lambda^2||x||^2$$.

Therefore in general we have that
$$||ax||=\lambda ||ax||$$.

If $0\neq x$ is an alternative element then $||yx||=||y||||x||$ for all $y$ (th 3.2).Thus if $x\in V_\lambda$, then by (1), $\lambda =1$,
and we are done with (2).

\hfill Q.E.D.

{\bf Remark} There are elements $0\neq a\in \widetilde{\aa}_n$ with $1\in$ Spec$(a)$ and $a$ non--alternative (See example above).
\vglue.15cm

Moreover now we prove that for each nonzero $a \in \widetilde{\aa}_n$,
$$
1\in \;\;\hbox{\rm Spec}\;\; (a)
$$
For this we need:
\vglue .5cm
\noindent
{\bf Lemma 3.4} For $a\neq 0$ in $\widetilde{\aa}_n\;\;n\geq 4$ there exists $x\neq 0$ in $\ache^\perp_a$ such that $||ax||=||a||||x||$.
\vglue .25cm
\noindent
{\bf Proof:} Recall that if $0\neq y\in \aa_n$ is alternative then $||ay||=||a||||y||$. (th 3.2).Thus
If $0\neq y \in \ache^\perp_a$ for some $y$ alternative we are done.

On the other hand by dimensional reasons (recall dim $\ache_a=4$ and any canonical basic element is alternative) no all alternative element belongs to $\ache_a$ so there exists $0\neq y\in \aa_n$ alternative element with $y=(y_1+y_2)\in \ache_a\oplus \ache^\perp_a$ and $y_2\neq 0$ in $\ache^\perp_a$.

But $||ay||=||a||||y||$ implies that
$$
\langle ay, ay\rangle=\langle a,a\rangle \langle y , y\rangle
$$
But 
\begin{eqnarray*}
\langle ay, ay\rangle&=&\langle a(y_1+y_2), a(y_1+y_2)\rangle\\
&=& \langle ay_1+ay_2, ay_1+ay_2\rangle\\
&=&\langle ay_1, ay_1\rangle +\langle ay_2, ay_2\rangle +2\langle ay_1, ay_2\rangle
\end{eqnarray*}
But $\langle ay_1, ay_2\rangle =-\langle a(ay_1), y_2\rangle =-\langle a^2y_1, y_2\rangle =0$ and $\langle ay_1, ay_1\rangle =\langle a, a\rangle \langle y_1, y_1\rangle =||a||^2||y_1||^2$ because $y_1\in \ache_a$  wich is associative and $y_2\perp y_1$; therefore
\begin{eqnarray*}
\langle ay, ay\rangle &=&\langle a, a\rangle \langle y_1, y_1\rangle +\langle ay_2, ay_2\rangle \qquad\hbox{\rm and}\\
\langle a, a\rangle \langle y, y\rangle &=&\langle a, a\rangle \langle y_1+y_2, y_1+y_2\rangle\\
&=&\langle a,a\rangle \langle y_1, y_1\rangle +\langle a, a\rangle \langle y_2, y_2\rangle
\end{eqnarray*}
because $y_1\perp y_2$.

Cancelling out we have $\langle a,a \rangle \langle y_2, y_2\rangle =\langle ay_2, ay_2\rangle$ and $||a||||y_2||=||ay_2||$ for $y_2\neq 0$ in $\ache^\perp_a$

\hfill Q.E.D.
\vglue .25cm
\noindent
{\bf Remark:} We can also use the construction of bilinear maps as in Adem (see [A],[L] and [Sh]) to prove this lemma because any non--zero double pure element belongs to a alternative subalgebra of $\aa_n$ with normed product (see [M2]).
\vglue .5cm
\noindent
{\bf Theorem 3.5} For $0\neq a\in \aa_n\;\;n\geq 3$ with $a$ doubly pure
$$
1\in\;\hbox{\rm Spec}(a).
$$
{\bf Proof:} Suppose that $1\not\in$ Spec$(a)$ then $a(ax)\neq a^2x$ for all (non-zero)  $x\in \ache^\perp_a$ then $a(ax)-a^2x\neq 0$ for all (non-zero) $x\in \ache^\perp_a$. Now $L^2_a(\ache_a)\subset \ache_a$ and $L^2_a(\ache^\perp_a)\subset \ache^\perp_a$ so we have that $L^2_a|_{\ache ^\perp_a}$ and  $(L^2_a-a^2I)|_{\ache ^\perp_a}$ are symmetric linear transformations.

Therefore $\langle a(ax)-a^2x,x\rangle \neq 0$ for all non-zero $x\in \ache^\perp_a$ and $\langle ax, ax\rangle\neq |a|^2\langle x, x\rangle$ and $||ax||\not =||a||||x||$ for all non-zero $x\in \ache^\perp_a$ which contradicts the previous lemma 3.4.

\hfill Q.E.D.
\noindent
\vglue .5cm

{\bf Corollary 3.6} For $a$ nonzero element in $\a_n$ $n\geq 3$ we have that:
$$dim Ker(a^2I - L^2_a)\geq 8.$$
\vglue .5cm
{\bf Proof}. Since any associator vanish when one of the entries is real,the question can be reduced to the case $a$ is a pure element.

On the other hand by Corollary 1.5  in [M2] is known that for $b$ double pure element
$$(b,\widetilde{e}_0,x)+(\widetilde{e}_0,b,x)=0$$

for all $x$ in $\a_n$.

Recalling that $\widetilde{e}_0$ is an alternative element, we can reduce the question to tha case $a$
double pure element,so by Theorem 3.5. $dim V_\lambda\geq 4$ for $\lambda=1$ and
$$dim(\ache_a\oplus V_1)\geq 8$$

\hfill Q.E.D.

\newpage
\noindent
{\Large\bf $\S$ 4. The zero divisors in $\aa_n\;\;\;n\geq 4$.}
\vglue.5cm

\vglue .25cm
In this section $0\neq a=(a_1, a_2)\in \aa_{n-1}\times \aa_{n-1}=\aa_n$ and $n\geq 4$\vglue .5cm
\noindent
{\bf Definition:} $a$ {\it is zero divisor if} there exists $x\neq 0$ in $\aa_n$ such that $ax=0$ i.e. 
$$
\hbox{\rm Ker}\;\; L_a\neq \{0\}.
$$
Since $||ax||=||xa||$ for all $x$ in $\aa_n$ (see [M1]) we have that $ax=0$ if and only if $xa=0$ so left and right zero divisors concide.

Also $||ax||=||\overline{a}x||$ (see [M1]) then if $x\neq 0\;\;t_n(a)x=(a+\overline{a})x=ax+\overline{a}x=0$ if $ax=0$; but $t_n(a)$ is real number so $ax=0 \Rightarrow t_n(a)=0$ so 
any zero divisor is a pure element.

Actually

{\bf Lemma 4.1} Any zero divisor in $\aa_n$ is double  pure.
\vglue .25cm
\noindent
{\bf Proof:} By Corollary 2.4 for $a\in \,_0\!\aa_n$
$$
\hbox{\rm Ker}\;L_a=\;\hbox{\rm Ker }L^2_a=\;\;\hbox{\rm Ker}\;L^2_{\widetilde{a}}=\;\;\hbox{\rm Ker}\;\;L_{\widetilde{a}}
$$
If $\{0\}\neq$ Ker$L_a$ then $\{0\}\neq$ Ker $L_{\widetilde{a}}$ and $\widetilde{a}=(-a_2, a_1)$ is a zero divisor so $\widetilde{a}$ is pure in $\aa_n$ and $a_2$ is pure in $\aa_{n-1}$ then $a$ is doubly pure.

\hfill Q.E.D.

Notice that $ax=0$ if and only if $||a||^{-1}ax=0$ so being a zero divisor is independent of the norm. On the other hand by Corollay 2.4(1) if $a$ is a zero divisor and  $r$ and $s$ in $\erre$ with $r^2+s^2=1$ then $L^2_a=L^2_{(a_1, a_2)}=L^2_{(ra_1 -sa_2, sa_1+ra_2)}=L^2_{ra+s\widetilde{a}}$ then $(ra+s\widetilde{a})$ is a zero divisor (in particular $\widetilde{a}$ is a zero divisor) and for any $b\neq 0$ in $\ache_a$ doubly pure we have that $b$ is a zero divisor and Ker$L_a=$ Ker$L_b$. Because 
\begin{eqnarray*}
L^2_{ra+s\widetilde{a}}&=&r^2L^2_a+s^2L_{\widetilde{a}}+rs(L_aL_{\widetilde{a}}+L_{\widetilde{a}}L_a)\\
&=&(r^2+s^2)L^2_{\widetilde{a}}+0\qquad\hbox{\rm (Corollary 2.5)}\\
&=&L^2_a.
\end{eqnarray*}
So if $||a||=1$ then $b=ra+s\widetilde{a}$ with $r^2+s^2=1$ and Ker $L_a=$ Ker$L_b$.
\vglue .5cm
\noindent
{\bf Lemma 4.2} Let $a=(a_1, a_2)$ be in $\aa_{n-1}\times \aa_{n-1}=\aa_n$ with $a\neq 0\;\;.a$ is a zero divisor if and only if $\hat{a}:=(a_2, a_1)$ is a zero divisor. Moreover $ax=0$ if and only if $\hat{a}\hat{x}=0$.
\vglue .25cm
\noindent
{\bf Proof:} Put $x=(x_1, x_2)$ so $\hat{x}=(x_2, x_1)$ in $_0\!\aa_{n-1}\times\,_0\!\aa_{n-1}$ calculating we have:
\begin{eqnarray*}
ax&=&(a_1, a_2)(x_1, x_2)=(a_1x_1+x_2a_2, x_2a_1-a_2x_1):=(c,d)\\
\hat{a}\hat{x}&=&(a_2, a_1)(x_2,x_1)=(a_2x_2+x_1a_1, x_1a_2-a_1x_2)=(\overline{c}, -\overline{d})
\end{eqnarray*}
therefore $c=d=0$ in $\aa_{n-1}$ if and only if $ax=\hat{a}\hat{x}=0$.

\hfill Q.E.D.
\vglue .25cm
\noindent
{\bf Remark:} Notice that lemma 4.1 said Ker$L_a=$ Ker$L_{\widetilde{a}}$ (equality) while lemma 4.2 said that Ker$L_a\cong$ Ker$L_{\hat{a}}$ (isomorphism), induced by 
$$_0\!\aa_{n-1}\times _0\!\aa_{n-1}\stackrel{\wedge}{\rightarrow}\,_0\!\aa_n\times \,_0\!\aa_n$$ 
given by $(x_1, x_2)\mapsto (x_2, x_1)$ which is a one--to one correspondence.

Now $0\neq a\in \aa_n$ is a zero divisor if and only if $ra$ is a zero divisor for all $r\in \erre$, $r\neq 0$ so we restrict ourselves to the study of zero divisor of norm one.

Using lemma 1.1 in $\S1$ we can construct zero divisors in the following way:
\vglue .5cm
\noindent
{\bf Proposition 4.3.} For $0\neq a\in \aa_n$, $n\geq 3$, and $a$ doubly pure we have:
\begin{enumerate}
\item[1)] If $b\in \ache^\perp_a$ and $|a|=|b|\neq 0$ then $(a,b)$ is a zero divisor in $\aa_{n+1}$.
\item[2)] $(a, \widetilde{a})$ and $(a, -\widetilde{a})$ are zero divisors in $\aa_{n+1}$.
\end{enumerate}
\noindent
{\bf Proof:}
\begin{enumerate}
\item[1)] For $b\in \ache^\perp_a$ with $|a|=|b|\neq 0$ we have that $(a,b)(\widetilde{a}, \widetilde{b})=(a\widetilde{a}-\overline{\widetilde{b}}b, \widetilde{b}a+b\overline{\widetilde{a}})=(a\widetilde{a}+\widetilde{b}b, \widetilde{b}a-b\widetilde{a})$ in $\aa_{n+1}$. By lemma 1.1. (2) $a\widetilde{a}+\widetilde{b}b=-||a||^2\widetilde{e}_0+||b||^2\widetilde{e}_0=0$ in $\aa_n$. By lemma 1.1. (6) $\;\;\widetilde{b}a=b\widetilde{a}\;\;$ because $b\in \ache^\perp_a$ in $\aa_n$.
\item[2)] Since $n\geq 3,\;\;\ache^\perp_a\not\neq\{0\}$. Taking $0\neq x \in \ache^\perp_a\subset \aa_n$
$$
(a, \widetilde{a})(x, -\widetilde{x})=(ax-\widetilde{x}\widetilde{a}, -\widetilde{x}a-\widetilde{a}x) \;\;\hbox{\rm in}\;\; \aa_{n+1}.
$$ 
By lemma 1.1. (5) $ax=\widetilde{x}\widetilde{a}$ and by lemma 1.1 (3) $-\widetilde{x}a-\widetilde{a}x=\widetilde{xa}+\widetilde{ax}=\widetilde{xa+ax}=0$ because $a\perp x$.

Similarly
$$
(a, -\widetilde{a})(x, \widetilde{x})=(ax-\widetilde{x}\widetilde{a}, \widetilde{x}a+\widetilde{a}x)=(0,0)\;\;\;\hbox{\rm in}\;\;\aa_{n+1}
$$
\hfill Q.E.D.
\end{enumerate}
\vglue .5cm
\noindent
{\bf Remark:} Notice that this proposition is telling us that any non--zero {\it doubly pure} element in $\aa_n\;\;(n\geq 3)$ is the component of a zero divisor in $\aa_{n+1}$.

We will prove that this is true also for non--zero {\it pure} elements in $\aa_n$.
\vglue .5cm
\noindent
{\bf Theorem 4.4.} For $0\neq a\in \aa_n$, $(n\geq 3)$ a doubly pure element of norm one and $0\neq \lambda \in$ Spec$(a)$ we have that $(a, \pm \lambda \widetilde{e}_0)$ are zero divisors in $\aa_{n+1}$.
\vglue .25cm
\noindent
{\bf Proof:} By definition of the spectrum of $a$
$$
a(ax)=-\lambda^2x\qquad\hbox{\rm for some}\qquad 0\neq x\in \ache^\perp_a
$$
Then in $\aa_{n+1}$
\begin{eqnarray*}
(a, \lambda \widetilde{e}_0)(ax, -\lambda \widetilde{x})&=&(a(ax)-\lambda^2\widetilde{x}\widetilde{e}_0, -\lambda \widetilde{x}a+\lambda \widetilde{e}_0 (\overline{ax}))\\
&=&(-\lambda ^2x-\lambda^2\widetilde{\widetilde{x}},\lambda (\widetilde{xa})-\lambda (\widetilde{xa}))\\
&=&(0,0)
\end{eqnarray*}
Similarly
\begin{eqnarray*}
(a, -\lambda \widetilde{e}_0)(ax, \lambda \widetilde{x})&=&(a(ax)+\lambda^2\widetilde{x}\widetilde{e}_0, \lambda \widetilde{x}a-\lambda \widetilde{e}_0(\overline{ax}))\\
&=&(-\lambda^2x-\lambda^2\widetilde{\widetilde{x}}, -\lambda \widetilde{xa}+\lambda (\widetilde{xa})\\
&=&(0,0).
\end{eqnarray*}
Notice also that $(a, \pm \lambda \widetilde{e}_0)(-1/\lambda ax, \pm \widetilde{x})=(0,0)$.

\hfill Q.E.D.

From this we derive the important:
\vglue .5cm
\noindent
{\bf Corollary 4.5} For $0\neq \alpha \in \aa_n\,(n\geq 3)$ with $|\alpha |=1$ and $\alpha$ pure element there exists $\beta$ in $\aa_n$ with $\alpha \perp \beta$ and $|\alpha|=|\beta|=1$ such that $(\alpha , \beta )$ is a zero divisor in $\aa_{n+1}$.
\vglue .5cm
\noindent
{\bf Proof:} If $\alpha$ is doubly pure we are done by Proposition 4.3. Suppose that $\alpha$ is pure element (no double pure) then there exist $a\in \aa_n$ doubly pure element and $r$ and $s$ real numbers with $r^2+s^2=1$ and $s\neq 0$ such that $\alpha =ra+s\widetilde{e}_0$ and $|a|=|\alpha |=1$.

By Corollary 2.4 (1) $L^2_{(a,-\widetilde{e}_0)}=L^2_{(ra+s\widetilde{e}_0, sa-r\widetilde{e}_0)}$. Therefore making $\beta =sa-r\widetilde{e}_0$ we have that
$$
L^2_{(a, -\widetilde{e}_0)}=L^2_{(\alpha , \beta)}
$$
Now by theorem 3.5 $\lambda =1\in$ Spec$(a)$ and by theorem 4.4 $(a, -\widetilde{e}_0)$ is a zero divisor in $\aa_{n+1}$ so $(\alpha , \beta)$ is a zero divisor in $\aa_{n+1}$.

By construction $|\alpha|=|\beta|=1$ and $\alpha \perp \beta$.

\hfill Q.E.D.
\vglue .5cm
{\bf Example:} We know that 2 is in  $Spec((e_1,e_2))$ in $\a_4$ so if $a=(\sqrt{2})^{-1}(e_1,e_2)$
and $\widetilde{e}_0=e_8$ in $\a_4$,according with Theorem 4.4 we have that
$$(a,2e_8)=\sqrt{2}^{-1}(e_1 +e_{10})+2e_{16}$$

is a zero divisor in $\a_5$.
\vglue.5cm

{\bf Remarks:} Notice that, in contrast, with the case $n=4$ where the zero divisors must have coordinates in $\a_3$ of equal norm (see [M1]) this is no the case for $\a_5$.(Besides the "trivial" cases when one of the coordinates is a zero divisor  and the other is equal to $0$ in $\a_4.$)

Therefore the zero divisors in $\a_5$ are "very far" to be described as in $\a_4$ where they can be identified with $V_{7,2}$ the Stiefel Manifold of 2 frames in $\erre^7.$

But also, Corollary 4.5 is telling us that the set of zero divisors in $\a_{n+1}$ has some subset wich 
can be describe in terms of the Stiefel Manifold $V_{2^n -1,2}.$ for $n\geq 3.$

\newpage
\noindent
{\Large\bf $\S$5.Zero divisors and Stiefiel manifolds.}
\vglue .5cm
\noindent
{\bf Definition:} $\alpha =(a,b)\in \,_0\!\aa_n \times \,_0\!\aa_n=\widetilde{\aa}_{n+1}$ is a {\it Stiefel element} if $||a||=||b||\neq 0$ and $a\perp b$ in $_0\!\aa_n$. 
\vglue .25cm
\noindent
{\bf Definition:} $\alpha =(a,b)\in \widetilde{\aa}_n\times \widetilde{\aa}_n\subset \widetilde{\aa}_{n+1}$ is a {\it Non--trivial element}  if $||a||=||b||\neq 0$ and $b\in \ache^\perp_a$ i.e. $a\perp b$ and $\widetilde{a}\perp b$.

 Clearly any non--trivial element in $\widetilde{\aa}_{n+1}$ is a Stiefel element in $\widetilde{\aa}_{n+1}$.
\vglue .5cm
\noindent
{\bf Lemma 5.1.} If $\alpha =(a,b)$ and $\chi =(x,y)$ are in $_0\!\aa_n \times \,_0\!\aa_n =\widetilde{\aa}_{n+1}$ then 
$$\langle \alpha, \chi \rangle=\langle a, x\rangle +\langle b, y\rangle$$
\vglue .25cm
\noindent
{\bf Proof:} 
\begin{eqnarray*}
\alpha \chi + \chi \alpha&=&(ax+yb, ya-bx)+(xa+by, bx-ya)\\
&=& (ax+xa+by+yb, 0)
\end{eqnarray*}
so
$$
-2\langle \alpha , \chi \rangle =\alpha \chi  + \chi \alpha =(-2\langle a, x\rangle -2 \langle y,b\rangle , 0)
$$
\hfill Q.E.D.
\vglue .5cm
\noindent
{\bf Theorem 5.2} $\alpha =(a,b)$ is a Stiefel element in $\widetilde{\aa}_{n+1}$ if and only if $(\alpha , \hat{\alpha})$ is a non--trivial element in $\widetilde{\aa}_{n+2}$.
\vglue .25cm
\noindent
{\bf Proof:} Recall that if $\alpha =(a,b)$ then $\hat{\alpha}=(b,a)$. By lemma 5.1 $\langle \alpha, \hat{\alpha}\rangle=\langle a, b\rangle +\langle b, a\rangle $ so $\alpha \perp \hat{\alpha}$ in $\widetilde{\aa}_{n+1}$ if and only if $a\perp b$ in $_0\!\aa_n$.

Also by lemma 5.1 $\langle \widetilde{\alpha}, \hat{\alpha}\rangle=\langle (-b,a), (b,a)\rangle =-||b||^2+||a||^2$ so $\widetilde{\alpha} \perp \hat{\alpha}$ in $\widetilde{\aa}_{n+1}$ if and only if $||a||=||b||\neq 0$.

\hfill Q.E.D.
\vglue .5cm
\noindent
{\bf Theorem 5.3} The set of non--trivial elements in $\widetilde{\aa}_{n+1}$ with entries of norm $1$ is the complex Stiefel manifold $W_{m,2}$ for $m=2^{n-1}-1$.[J].
\vglue .25cm
\noindent
{\bf Proof:} Define ${\cal H}_n:\widetilde{\aa}_n \times \widetilde{\aa}_n\rightarrow \cee =$ complex numbers, by ${\cal H}_n (x,y)=2\langle x, y\rangle -2i\langle \widetilde{x}, y\rangle$.
\vglue .25cm
\noindent
{\bf Claim:} ${\cal H}_n$ define a Hermitian inner product in $\widetilde{\aa}_n$.

 Clearly ${\cal H}_n$ is $\erre$--bilinear and
\begin{eqnarray*}
\overline{{\cal H}_n (x,y)} &=& 2\langle x,y\rangle +2i \langle \widetilde{x}, y\rangle =2\langle x , y \rangle -2i\langle \widetilde{y}, x\rangle\\
&=&{\cal H}_n (y,x).
\end{eqnarray*}
Also
\begin{eqnarray*}
{\cal H}_n (\widetilde{x}, y) &=&2\langle \widetilde{x}, y\rangle -2i \langle \widetilde{\widetilde{x}}, y\rangle =2\langle \widetilde{x}, y\rangle +2i \langle x,y\rangle\\
&=& 2i \langle x, y \rangle -2i^2\langle \widetilde{x}, y\rangle =i {\cal H}_n (x,y)
\end{eqnarray*}
Therefore ${\cal H}_n$ define a Hermitian product as claimed. 

By definition ${\cal H}_n (a,b)=0$ if and only if $b\in \ache^\perp_a$ and 
$$
W_{m,2}=\{(a,b)\in \widetilde{\aa}_n\times \widetilde{\aa}_n|{\cal H}_n (a,b)=0\;\;\hbox{\rm and}\;\;||a||=||b||=1\}
$$
for
$$
m={1\over 2} (2^n-2)=2^{n-1}-1
$$
\hfill Q.E.D.
\vglue .25cm
\noindent
{\bf Remark} Notice that the set of Stiefel elements in $\widetilde{\aa}_{n+1}$, with entries of norm one, in $_0\!\aa_n$ can be seen as the real Stiefel manifold $V_{2^n-1,2}$ i.e. 
$$
V_{2^n-1,2}=\{(a,b)\in \,_0\!\aa_n\times \,_0\!\aa_n |a\perp b\;\;\;\hbox{\rm and}\;\;|a|=|b|=1\}
$$

Now we give a partial answer to the following:
\vglue.25cm
\noindent
{\bf Question:} Is any Stiefel element in $\widetilde{\aa}_{n+1}$ a zero divisor?.
\vglue .25cm
\noindent
{\bf Partial Answer:} So far, we have seen that the following types of Stiefel elements are zero divisors.

Suppose that $\alpha =(a,b)\in \widetilde{\aa}_{n+1}$ with $a\in \, _0\aa_n, b\in \, _0\aa_n$ and $|a|=|b|\neq 0$ with $a\perp b$.
\vglue .25cm
\noindent
{\bf Case: $a$ and $b$ doubly pure} in $\aa_n$.
\begin{enumerate}
\item If $\alpha =(a,b)\in \widetilde{\aa}_n\times \widetilde{\aa}_n$ is a non--trivial then $b\in \ache^\perp_a$ and by $\S$4 
$(a,b)(\widetilde{a}, \widetilde{b})=(0,0)$ in $\aa_{n+1}.$

\item Suppose that $b\in \ache_a$ .

Since $|a|=|b|\neq 0\;\;\;b=ra+s\widetilde{a}$ for $r^2+s^2=1$, and $a\perp b$ implies that $0=\langle a, b\rangle =\langle a, ra\rangle +\langle a, s\widetilde{a}\rangle$,but $a\perp \widetilde{a}$ and  we have that $0=r\langle a, a\rangle =r||a||^2$, therefore $r=0$ and $b=s\widetilde{a}$ with $s^2=1$ i.e. $b=\pm\widetilde{a}$.
\end{enumerate}
In $\S$ 4 we show that $(a, \pm \widetilde{a})(x, \mp \widetilde{x})=(0,0)$ in $\aa_{n+1}$ for all $0\neq x\in \ache^\perp_a$, so $(a,b)$ is a zero divisor.
\vglue .25cm
\noindent
{\bf General case:} $a$ and $b$ pure elements in $\aa_n$.

In $\S$ 4 we prove that for any non--zero doubly pure element $c$ in $\widetilde{\aa}_n\;\;\;(c, \widetilde{e}_0)$ is a zero divisor in $\aa_{n+1}$.

Since any pure element $a\in \, _0\aa_n$ with $a\neq 0$ is of the form $a=rc\mp s\widetilde{e}_0$ for $r,s$ in $\erre$ with $r^2+s^2=1$ and $|a|=|c|\neq 0$ for some $c$ doubly pure in $\widetilde{\aa}_n$ then $(rc\mp s\widetilde{e}_0, sc\pm r\widetilde{e}_0)$ is a zero divisor so $(a,b)$ is a zero divisor in $\aa_{n+1}$ for $a\in \, _0\aa_n\quad a\neq 0$ and $b=sc\pm r\widetilde{e}_0$ if $a=rc\mp s\widetilde{e}_0$.

More generally if $0\neq\lambda\in Spec(c)$ then $(c,\lambda\widetilde{e}_0)$ is a zero divisor in $\aa_{n+1}$ for $c$ a doubly pure element in $\widetilde{\aa}_n$ with $c\neq 0$ (\S 4). If $a=rc\mp s\lambda\widetilde{e}_0$ with $r$ and $s$ in $\erre$ such that
$$r^2+(s\lambda)^2=1$$
then $(a,b)$ is a zero divisor for $b=(s\lambda)c\pm r\widetilde{e}_0$. 
\vglue .1cm
{\bf Open Question:}

Given $a\in \,_0\aa_n$ with $a\neq 0$.

If $b\in \, _0\aa_n$ with $b\perp a, |b|=|a|\neq 0$ then 

Is $b$ of one the forms described above?

{\bf Remark.}Notice that if $\lambda=0$ is in $Spec(c)$ then $c$ is a zero divisor and $a=rc\mp s\lambda\widetilde{e}_0=rc$ with $r^2=1$ so $a=\pm c$ and $b=\mp \widetilde{e}_0$ and $(a,b)$ is a zero divisor.
\vglue 1cm
\centerline{\bf REFERENCES}
\begin{enumerate}
\item[{\rm [A]}] J. Adem, {\it Construction of some normed maps}, Bol. Soc.Mat. Mex. 20 (1975) 59--75.
\item[{\rm [B-D-I]}] Biss -Dugger-Isaksen,{\it Large annihilators in Cayley-Dickson algebras}
RA/0511691 Nov 2005.
\item[{\rm [Ch]}] Chaitin-Chatelin F,{Inductive multiplication in Dickson algebras} CERFACS  Technical
Report TR/PA/05/56.
\item[{\rm [E-S]}] Eakin--Sathaye, {\it On automorphisms and derivations of Cayley--Dickson algebras}, Jornal of Pure and Applied Algebra, 129, 263--278 (1990).
\item[{\rm [H-J]}] Horn and Jhonson, {\it Matrix analysis}, Cambridge University Press.
\item[{\rm [J]}] James I., {\it Topology of Stiefel manifolds}, Cambridge University Press.
\item[{\rm [L]}] Lam, K.Y., {\it Topological methods for studying the composition of quadratic forms}, Canadian Math. Soc. Conf. Proc. 4 173--192 (1984).
\item[{\rm [M1]}] Moreno G., {\it The zero divisor of the Cayley--Dickson algebras over the real numbers}, Bol. Soc. Mat. Mex. (tercera serie), Vol. 4 No. 1 1998.
\item[{\rm [M2]}] Moreno G., {\it Alternative elements in the Cayley--Dickson algebras}.Proceedings in Relativity and Field theory in Honor of Jerzy Plebansky Birkhauser 2004.
math.RA/0404395.
\item[{\rm [M3]}] Moreno G. {\it Hopf construction in higher dimensions} Bol. Soc. Mat. Mex 2004
also math.AT/0404172.
\item[{\rm [M4]}] Moreno G. {\it Monomorphisms between Cayley--Dickson algebras}.Nonassocitive 
algebras and its applications Chapter 22 CRC Press, Taylor and Francis 2005.math.RA/0512516.
\item[{\rm [Po]}] Porteous, I., {\it Topological algebra}, Cambridge University Press.
\item[{\rm [Pr]}] Praslov, {\it Problems and theorems in linear algebra.} Translation of the A.M.S., Vol. 134. 
\item[{\rm [Sch1]}] Schaffer, R.D., {\it An introduction to non--associtive algebras}, Academic Press 1966.
\item[{\rm [Sch2]}] Schaffer, R.D., {\it On algebras formed by the Cayley Dickson Process}, American J. of Math. 76 (1954) 435--446.
\item[{\rm [Sh]}] Shapiro, D., {\it Product of sum of squares}, Exposition Math. (2) 235--261, (1984).
\item[{\rm [K-Y]}]Khalil-- Yiu P., {\it The Cayley--Dickson algebras: A theorem of Hurwitz and quaternions.}, Bulletin de la Societ\'e, de L\'odz, Vol. XLVIII (1997)pp. 117--169.
\item[{\rm [Z]}] Zhang F., {\it Matrix theory}, Springer--Verlag (University Text).
\end{enumerate}
\newpage

{\sc \begin{flushleft} 
   Guillermo Moreno Rodr\'{\i}guez \newline 
    Departamento de Matematicas \newline 
  Centro de Investigaci\'on y de estudios Avanzados      \newline 
   (CINVESTAV del IPN) Apartado Postal 14--740             \newline 
   MEXICO  D.F. MEXICO 07300                         \newline 
   {\rm gmoreno@math.cinvestav.mx}

         \end{flushleft} } 

\end{document}